\documentclass[12pt,letterpaper]{article}
\usepackage{amssymb, amsthm, amsmath}
\usepackage{hyperref}
\usepackage[margin=1in]{geometry}

\usepackage[boxruled]{algorithm2e}

\usepackage{tikz}
\tikzset{every picture/.style={thick}}
\tikzset{every node/.style={circle,draw=black,inner sep=2pt}}
\usepackage{forest}

\def\FIRSTEXAMPLE{\begin{tikzpicture}[scale = .5]
	
	\node (1) at (0,2) {1};
	\node (2) at (0,0) {2};
	\node (3) at (0,-2) {3};
	\node (4) at (2,0) {4};
	\node (5) at (4,0) {5};
	\node (6) at (6,2) {6};
	\node (7) at (6,0) {7};
	\node (8) at (6,-2) {8};
	
	\draw (4) -- (1);
	\draw (4) -- (2);
	\draw (4) -- (3);
	\draw (4) -- (5);
	\draw (5) -- (6);
	\draw (5) -- (7);
	\draw (5) -- (8);

\end{tikzpicture}}

\def\KBOUND{\begin{tikzpicture}
\node (0) at (0,0) {};
\foreach \i in {1,...,4}{
\pgfmathsetmacro{\ang}{45+90*(\i-1)}
\pgfmathsetmacro{\sang}{\ang-20}
\pgfmathsetmacro{\lang}{\ang+20}
\pgfmathsetmacro{\sind}{int(4+2*(\i-1))}
\pgfmathsetmacro{\lind}{int(4+2*(\i-1)+1)}
\node (\i) at (\ang:1) {};
\draw (\i) -- (0);
\node (\sind) at (\sang:1.4) {};
\node (\lind) at (\lang:1.4) {};
\draw (\sind) -- (\i) -- (\lind);
}
\end{tikzpicture}}

\def\COMB{\begin{tikzpicture}
\node (0) at (2,2) {};
\node (1) at (2,1) {};
\foreach \x in {2,...,8}{
\node (\x) at (\x,0) {};
}
\foreach \x in {3,4,5,7}{
\node (\x-1) at (\x,1) {};
}
\foreach \x in {3,5,7}{
\node (\x-2) at (\x,2) {};
}
\foreach \x in {3}{
\node (\x-3) at (\x,3) {};
}
\draw (0) -- (1) -- (2) -- (3) -- (4) -- (5) -- (6) -- (7) -- (8);
\draw (3) -- (3-1) -- (3-2) -- (3-3);
\draw (4) -- (4-1);
\draw (5) -- (5-1) -- (5-2);
\draw (7) -- (7-1) -- (7-2);

\node[draw=none,rectangle] at (8.5,0) {$v$};
\node[draw=none,rectangle] at (1.5,0) {$u$};
\end{tikzpicture}}

\def\PICKCOMB{\begin{tikzpicture}
\foreach \x in {1,...,5}{
\node (\x-0) at (\x,0) {};
\node (\x-1) at (\x,1) {};
}
\draw (1-0) -- (2-0) -- (3-0) -- (4-0) -- (5-0);
\draw (1-1) -- (2-1) -- (3-1) -- (4-1) -- (5-1);
\draw (1-0) -- (3-1) -- (2-0) -- (4-1) -- (5-0);
\draw (1-1) -- (1-0) -- (2-1);
\draw (4-0) -- (4-1);

\begin{scope}[yshift=5.5cm,xshift=7cm]
\node (0) at (2,-2) {};
\node (1) at (1,-2) {};
\foreach \x in {2,...,8}{
\node (\x) at (0,-\x) {};
}
\foreach \x in {3,4,5,7}{
\node (c\x-1) at (1,-\x) {};
}
\foreach \x in {3,5,7}{
\node (c\x-2) at (2,-\x) {};
}
\foreach \x in {3}{
\node (c\x-3) at (3,-\x) {};
}
\draw (0) -- (1) -- (2) -- (3) -- (4) -- (5) -- (6) -- (7) -- (8);
\draw (3) -- (c3-1) -- (c3-2) -- (c3-3);
\draw (4) -- (c4-1);
\draw (5) -- (c5-1) -- (c5-2);
\draw (7) -- (c7-1) -- (c7-2);
\end{scope}

\draw (8) -- (5-0) -- (2) -- (5-1);
\draw (2) to[bend right] (8);

\node[draw=none,rectangle] at (6.5,3.5) {$u$};
\node[draw=none,rectangle] at (6.5,-2.5) {$v$};

\end{tikzpicture}}

\def\ALGB{\begin{tikzpicture}
\node[rectangle,inner sep=3pt] (0) at (-1.5,1) {\small$2$};
\node[rectangle,inner sep=3pt] (1) at (-2,0) {\small$2$};
\node[rectangle,inner sep=3pt] (2) at (-1,0) {\small$0$};
\node[rectangle,inner sep=3pt] (4) at (-2.8,-1) {\small$0$};
\node[rectangle,inner sep=3pt] (5) at (-2,-1) {\small$0$};
\node[rectangle,inner sep=3pt] (6) at (-1.2,-1) {\small$0$};
\node[rectangle,inner sep=3pt] (3) at (0.8,2.5) {\small$4$};
\node[rectangle,inner sep=3pt] (7) at (0.8,1) {\small$2$};
\node[rectangle,inner sep=3pt] (8) at (3.2,1) {\small$2$};
\node[rectangle,inner sep=3pt] (9) at (0,0) {\small$0$};
\node[rectangle,inner sep=3pt] (10) at (0.8,0) {\small$0$};
\node[rectangle,inner sep=3pt][rectangle,inner sep=3pt] (11) at (1.6,0) {\small$0$};
\node[rectangle,inner sep=3pt] (12) at (2.4,0) {\small$0$};
\node[rectangle,inner sep=3pt] (13) at (3.2,0) {\small$0$};
\node[rectangle,inner sep=3pt] (14) at (4,0) {\small$0$};

\draw[->] (3) to[out=-90,in=90] (0);
\draw[->] (3) to[out=-90,in=90] (7);
\draw[->] (3) to[out=-90,in=90] (8);
\draw[->] (0) to[out=-90,in=90] (1);
\draw[->,dashed] (0) to[out=-90,in=90] (2);
\draw[->] (1) to[out=-90,in=90] (4);
\draw[->] (1) to[out=-90,in=90] (5);
\draw[->] (1) to[out=-90,in=90] (6);
\draw[->] (7) to[out=-90,in=90] (9);
\draw[->] (7) to[out=-90,in=90] (10);
\draw[->] (7) to[out=-90,in=90] (11);
\draw[->] (8) to[out=-90,in=90] (12);
\draw[->] (8) to[out=-90,in=90] (13);
\draw[->] (8) to[out=-90,in=90] (14);
\end{tikzpicture}}

\def\ALGA{\begin{tikzpicture}
\node[rectangle,inner sep=3pt] (0) at (0,3.5) {\small$3$};
\node[rectangle,inner sep=3pt] (1) at (-2,2) {\small$2$};
\node[rectangle,inner sep=3pt] (2) at (0,2) {\small$0$};
\node[rectangle,inner sep=3pt] (3) at (2,2) {\small$3$};
\node[rectangle,inner sep=3pt] (4) at (-2.8,1) {\small$0$};
\node[rectangle,inner sep=3pt] (5) at (-2,1) {\small$0$};
\node[rectangle,inner sep=3pt] (6) at (-1.2,1) {\small$0$};
\node[rectangle,inner sep=3pt] (7) at (0.8,1) {\small$2$};
\node[rectangle,inner sep=3pt] (8) at (3.2,1) {\small$2$};
\node[rectangle,inner sep=3pt] (9) at (0,0) {\small$0$};
\node[rectangle,inner sep=3pt] (10) at (0.8,0) {\small$0$};
\node[rectangle,inner sep=3pt] (11) at (1.6,0) {\small$0$};
\node[rectangle,inner sep=3pt] (12) at (2.4,0) {\small$0$};
\node[rectangle,inner sep=3pt] (13) at (3.2,0) {\small$0$};
\node[rectangle,inner sep=3pt] (14) at (4,0) {\small$0$};

\draw[->, dashed] (0) to[out=-90,in=90] (1);
\draw[->, dashed] (0) to[out=-90,in=90] (2);
\draw[->] (0) to[out=-90,in=90] (3);
\draw[->, dashed] (1) to[out=-90,in=90] (4);
\draw[->, dashed] (1) to[out=-90,in=90] (5);
\draw[->, dashed] (1) to[out=-90,in=90] (6);
\draw[->] (3) to[out=-90,in=90] (7);
\draw[->] (3) to[out=-90,in=90] (8);
\draw[->] (7) to[out=-90,in=90] (9);
\draw[->] (7) to[out=-90,in=90] (10);
\draw[->] (7) to[out=-90,in=90] (11);
\draw[->] (8) to[out=-90,in=90] (12);
\draw[->] (8) to[out=-90,in=90] (13);
\draw[->] (8) to[out=-90,in=90] (14);
\end{tikzpicture}}

\DeclareMathOperator{\nullity}{nullity}

\allowdisplaybreaks

\newtheorem{theorem}{Theorem}
\newtheorem{proposition}{Proposition}
\newtheorem{lemma}{Lemma}
\newtheorem{corollary}{Corollary}
\newtheorem{observation}{Observation}

\theoremstyle{definition}

\newtheorem{definition}{Definition}
\newtheorem{example}{Example}

\title{Properties of a $q$-analogue of zero forcing}
\author{
Steve Butler \and
Craig Erickson \and 
Shaun Fallat \and 
H. Tracy Hall \and 
Brenda Kroschel \and 
Jephian C.-H. Lin \and 
Bryan Shader \and 
Nathan Warnberg \and 
Boting Yang }
\date{\empty}

\begin{document}
\maketitle

\begin{abstract}
Zero forcing is a combinatorial game played on a graph where the goal is to start with all vertices unfilled and to change them to filled at minimal cost.  In the original variation of the game there were two options. Namely, to fill any one single vertex at the cost of a single token; or if any currently filled vertex has a unique non-filled neighbor, then the neighbor is filled for free.  This paper investigates a $q$-analogue of zero forcing which introduces a third option involving an oracle. Basic properties of this game are established including determining all graphs which have minimal cost $1$ or $2$ for all possible $q$, and finding the zero forcing number for all trees when $q=1$.
\end{abstract}

\section{Introduction}
The zero forcing game is a combinatorial game played on a graph.  The game involves filling in the vertices of a graph by certain legal moves, the most important of which is the following \emph{filling rule} (sometimes known as the forcing rule or coloring rule): If a filled vertex has a unique unfilled neighbor (and any number of filled neighbors), then the unfilled neighbor becomes filled.  The game is summarized as follows.

\begin{quote}
{\bf The Zero Forcing Game (or $Z$-Game)} -- All the vertices of the graph $G$ are initially unfilled and there is one player who has tokens.  The player will repeatedly apply one of the following two operations until all vertices are filled:
\begin{enumerate}
\item For one token, any vertex can be changed from unfilled to filled.
\item At no cost, the player can apply the filling rule.
\end{enumerate}
\end{quote}

The $Z$-Game number of a graph, denoted $Z(G)$, is the minimum number of tokens needed to guarantee that all vertices can be filled (this is sometimes referred to as the ``zero forcing number'' or ``fast-mixed search number'').  Zero forcing was developed in the combinatorial matrix theory community to give a bound for the minimum rank of a symmetric matrix associated with a graph \cite{AIM}, zero forcing has also been developed independently for other purposes (see \cite{FH} and references contained therein).  

\begin{theorem}[\cite{AIM}]
If $A$ is a real-symmetric matrix with nonzero off-diagonal entries corresponding to the edges of $G$, then $\nullity(A)\le Z(G)$.
\end{theorem}

There have been many variations on the $Z$-Game, the one which will be considered here is a $q$-analogue of zero forcing  which introduces a new operation available to the player (see \cite{BGH}).  This new operation allows the (potential) application of the filling rule on a \emph{smaller} part of the graph.  For $W\subseteq V(G)$, let $G[W]$ denote the induced subgraph of $G$ on the vertices $W$.

\begin{quote}
{\bf The $q$-Analogue of the Zero Forcing Game (or $Z_q$-Game)} -- All the vertices of the graph $G$ are initially unfilled and there is one player who has tokens, and one oracle.  The player will repeatedly apply one of the following three operations until all vertices are filled.
\begin{enumerate}
\item For one token, any vertex can be changed from unfilled to filled.
\item At no cost, the player can apply the filling rule.
\item Let the vertices currently filled be denoted by $F$, and $U_1,\ldots,U_k$ be the vertex sets of the connected components of $G[V\setminus F]$ (i.e., components of unfilled vertices).  If $k\ge q+1$, the player can select \emph{at least}\/ $q+1$ of the $U_i$ and announces the selection to the oracle.  The oracle selects a nonempty subset of these components, $\{U_{i_1},\ldots,U_{i_\ell}\}$, and announces it back to the player.  At no cost, the player can apply the filling rule on $G[F\cup U_{i_1}\cup\cdots\cup U_{i_\ell}]$. 
\end{enumerate}
\end{quote}

The $Z_q$-Game number of a graph, denoted $Z_q(G)$, is the minimum number of tokens needed to guarantee that all vertices can be filled, regardless of how the oracle responds.  Due to the ``at least'' part the following inequalities are obtained: $Z_0(G)\le Z_1(G)\le\cdots\le Z(G)$ (think of $Z(G)=Z_\infty(G)$ as never being able to appeal to the oracle).

The parameter $Z_q(G)$ also gives a bound related to maximum nullity.  

\begin{theorem}[\cite{BGH}]
If $A$ is a real-symmetric matrix with nonzero off-diagonal entries corresponding to the edges of $G$ \emph{and} $A$ has at most $q$ negative eigenvalues, then $\nullity(A)\le Z_q(G)$.
\end{theorem}

For the $Z$-Game it is known that the spending of tokens can all happen up front before applying the filling rule.  Hence, in the literature there is a focus on zero forcing \emph{sets}.  In the $Z_q$-Game it might be disadvantageous to spend all tokens up front, i.e., the oracle's response(s) may change the optimal spending pattern.

\begin{example}\label{ex:basic}
Consider the graph $T$ shown in Figure~\ref{fig:firstexample}.  For this graph $Z(T)=4$, e.g., tokens must be spent on at least two of $\{1,2,3\}$ and at least two of $\{6,7,8\}$ and, having done so, the filling rule can be used on the remaining vertices.

\begin{figure}[ht!]
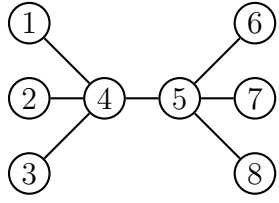

\centering
\FIRSTEXAMPLE
\caption{$Z(T) = 4$ and $Z_1(T) = 3$.}
\label{fig:firstexample}
\end{figure}

On the other hand, $Z_1(T)=3$.  To see this, initially spend on $1$ and $6$ and then apply the filling rule for $4$ and $5$.  At this point, hand the oracle $2$ and $7$.  Whatever is returned will be filled by the filling rule.  Continue handing one unfilled vertex from each side to the oracle until one side is completely filled.  At this point, at most $1$ more token must be spent to get down to one unfilled vertex on the other side (and if the oracle had consistently returned one side over the other this would be necessary).  Finally, the filling rule may be applied to fill in the last vertex.

Note if three tokens has been spent up front then there would exist a vertex with two adjacent unfilled leaves.  At this point the oracle could prevent these leaves from being filled.  So to achieve the optimal value, a delay in spending tokens is necessary.
\end{example}

In Example \ref{ex:basic}, the use of the response of the oracle allowed a reduction in spending on one side of the graph by one.\footnote{In terms of the linear algebraic philosophy of zero forcing, tokens were spent to probe and get additional information about the possible structure of the null space of a matrix and using the information to obtain a better bound on the nullity.}  In particular, there is a need to shift from zero forcing sets and into zero forcing strategies.

While there has been work on the combinatorial aspects of the $Z$-Game (see \cite{FH}), less is known about the $Z_q$-Game.  This paper begins to address the situation. Section~\ref{sec:small} considers what happens when $Z_q(G)$ is small compared to $q$.  Section~\ref{sec:PickComb} determines all connected graphs $G$ which have $Z_q(G)=1$ or $2$.  Section~\ref{sec:Z1onTree} gives an efficient method to compute $Z_1(G)$ on trees.  Finally, some concluding remarks and paths of future investigation are provided.

\section{Small forcing number} \label{sec:small}
In the $Z_q$-Game, in order for the player to use the option involving an oracle, there is a need to have sufficiently many components which are ``far'' from each other.  Intuitively, this means that a ``fair amount'' of tokens will already have been spent before the oracle operation is useful.  Lemma \ref{lem:spendk} makes this observation more precise.

\begin{lemma}\label{lem:spendk}
Consider the $Z_q$-Game on a graph $G$.  If at most $q$ tokens have been spent, there is no place on which the filling rule can be applied, and the oracle has not yet returned a set of components on which filling can occur, then for any group of at least $q+1$ components announced to the oracle there is a nonempty subset which the oracle can respond with and for which no filling occurs.
\end{lemma}

\begin{proof}
Let $F = \{f_1,f_2,\dots,f_k\}$ be the set of filled vertices that are adjacent to at least one unfilled vertex.  Based on the assumption that no forcing has occurred on components returned by the oracle it follows that $k\le q$.  This is because the number of such vertices equals the number of tokens spent.  (Note that it is not necessarily the case that the $f_i$ are where the tokens have been spent, as it is possible that the filling rule was applied.  However, any such vertex used to apply the filling rule cannot be one of the $f_i$.)

Now, let $\{\mathcal{C}_1,\mathcal{C}_2,\dots,\mathcal{C}_\ell\}$ with $\ell\ge q+1$ be the connected components of $G[V(G) \setminus F]$ that are handed to the oracle.  Consider the incidence array $A$ with rows corresponding to the vertices $f_i$, for $1\le i \le k$, and the columns corresponding to the connected components $\mathcal{C}_{j}$, for $1\le j \le \ell$, where the $i,j$ entry is 
\[
A_{i,j} = \left\{\begin{array}{c@{\qquad}l} 1 & \text{if $f_i$ is adjacent to $\mathcal{C}_{j}$,} \\ 0 & \text{else.} \end{array}\right.
\]
Here adjacent means there is some vertex in $C_j$ adjacent to vertex $f_i$.

The interest lies in finding a set of columns in this array with the property that if the corresponding components are returned then no filling rule can be used.  Note that any column with all $0$s can be returned as none of the $f_i$ are adjacent to the component and so no filling rule can be used.  So, without loss of generality, assume that the columns are all nonzero.

Repeat the following action as much as possible: If any row has row sum $1$ then delete that row and the corresponding column where the $1$ entry occurred.  Once this action is no longer possible, whichever columns remain are the components the oracle will return.

To see why, note that for every remaining row (e.g., vertex $f_i$) that it is adjacent to either $0$ of the components or adjacent to $2$ or more.  In either situation the filling rule does not apply.  Further, for any row which was deleted, the remaining columns all had to have $0$ entries as the only column with an entry of $1$ was removed; hence those $f_i$ are not adjacent to the components returned.

It remains to show that there is \emph{something} available to hand back.  For this, notice that $\ell\ge q+1>q\ge k$ so when the row/column elimination algorithm is repeated not every column will get deleted and if the array reduces down to one row there will be at least two columns left.  Further, if the array reduces to one row it must contain all ones.  This is because of the assumption that  no column had all zeroes so the remaining columns had at least one nonzero entry and the only entries which would have been removed to this point have value zero.
\end{proof}

By Lemma \ref{lem:spendk}, if $Z_q(G) \le q$, then there has to be a way to spend tokens without using the oracle as the oracle can always return a set of components on which no filling rule can be used.  As a consequence, it must be possible to fill all vertices by spending tokens and applying the filling rule, in other words, using only the methods allowed in the $Z$-Game.  This establishes Theorem \ref{thm:spendfirst}.

\begin{theorem}\label{thm:spendfirst}
If $Z_k(G)\le k$, then $Z_k(G)=Z(G)$.
\end{theorem}

\begin{corollary}
If $Z_k(G)=k$, then $Z_{k-1}(G)=k$.
\end{corollary}
\begin{proof}
By Theorem~\ref{thm:spendfirst}, $Z(G)=k$.  Also, $Z_{k-1}(G)\le Z_k(G)\le k$.  If $Z_{k-1}(G)\le k-1$, then $Z(G)=Z_{k-1}(G)\le k-1$ by Theorem~\ref{thm:spendfirst}
, which is impossible.
\end{proof}

The assumption in Theorem~\ref{thm:spendfirst} cannot be changed to $Z_k(g)\le k+1$ as shown by the following example.

\begin{proposition}\label{prop:kbound}
Let $T$ be the tree on $3k+4$ vertices formed by taking $k+1$ copies of $K_{1,3}$ and gluing them together on a leaf.  Then $Z(T)=k+2$ and $Z_k(T)=k+1$. 
\end{proposition}

\begin{figure}[ht!]
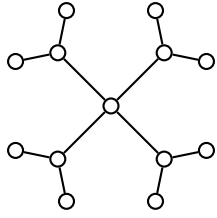

\begin{center}
\KBOUND
\end{center}
\caption{The graph from Proposition~\ref{prop:kbound} with $k+1=4$.}
\label{fig:kbound}
\end{figure}

\begin{proof}
The process begins by showing $Z(T) = k+2$ by establishing the path cover number (for the $Z$-Game on a tree the path cover number is the zero forcing number \cite{FH}). Since there are $2k+2$ leaves and each path covers at most $2$ leaves so the graph has to have at least $k+1$ paths.  If there are exactly $k+1$ paths, then either the center vertex is not covered by a path or a path covers the central vertex and then there is a leaf which is isolated and so more than $k+1$ paths would be needed to cover all leaves.  Thus, at least $k+2$ paths are needed to cover $T$.  On the other hand such a path cover is easily found, e.g., take $k+1$ paths of three vertices (each consisting of two pendant vertices and a degree $3$ vertex in the middle) and the path consisting of the center vertex.

By Theorem~\ref{thm:spendfirst}, $Z_k(T)\ge k+1$, since otherwise we would have $Z(T)\le k$.  It now suffices to give a strategy where $k+1$ tokens are spent and can fill all vertices.  To do this, spend $1$ token in a leaf in each copy of $K_{1,3}$.  Now, apply the filling rule to fill all the degree $3$ vertices.  At this point there are $k+1$ unfilled leaves which are handed to the oracle, whatever gets returned gets filled.  In particular, there is some filled degree $3$ vertex whose leaf neighbors are filled which, by the filling rule, fills the degree $k+1$ vertex.  Applying the filling rule again to all remaining leaves fills the entire graph.
\end{proof}

\section{Graphs with \texorpdfstring{$Z_q(G)$}{Zq(G)} equal to \texorpdfstring{$1$}{1} or \texorpdfstring{$2$}{2}}\label{sec:PickComb}

When $Z_q(G)=1$ there are two possibilities, either $q=0$ or $q\ge 1$ (in the latter case $Z_q(G)=Z(G)$.  When $q=0$, the oracle does not enter into the game as the filling rule can be applied to each component independently (this is known as the positive semi-definite $Z$-Game or the $Z_+$-Game, see \cite {EGR}).  In particular, if $G$ is a tree, then after spending one token the filling rule can be applied to each component and fill the entire graph.  On the other hand, if the graph contains a cycle, then the filling propagation will stop on a component which has a cycle.  So $Z_0(G)=1$ if and only if $G$ is a tree (see \cite{BBF}).  

When $q\ge 1$, if only one token is needed to fill the entire graph, then it must be the case that after the token is spent the filling rule works for the rest of the graph.  Following the propagation of the filling process, the structure must be a path (e.g., there is always at most one unfilled neighbor of each vertex).

This establishes Theorem \ref{thm:Z0}.
\begin{theorem}[\cite{BBF}]\label{thm:Z0}
Let $G$ be a connected graph.  Then,
\begin{itemize}
\item $Z_0(G)=1$ if and only if $G$ is a tree.
\item $Z_q(G)=1$ for $q\ge 1$ if and only if $G$ is a path.
\end{itemize}
\end{theorem}

To understand the case when $Z_q(G)=2$ it is necessary to understand how the oracle can respond in the $Z_1$-game to impede progress of the player.

\begin{definition}
A \emph{fort} in a graph $G$ is a subset $W\subseteq V(G)$ such that the induced subgraph $G[W]$ has at most two connected components and every vertex outside of $W$ has either zero or two or more neighbors in $W$.
\end{definition}

\begin{observation}
Let $W$ be a fort in a graph $G$.  If at some point in a $Z_1$-game every vertex in $W$ is unfilled, then the player cannot win without spending at least one more token.
\end{observation}

The filling rule cannot proceed in $W$ since each filled neighbor of a vertex in $W$ is adjacent to multiple vertices in $W$.  If the player gives at least two components to the oracle then the oracle either returns components containing nothing in $W$ or containing all of $W$.  In any case no filling in $W$ can occur so the player must spend at least one more token.

Note that this does not say that a token must be spent in \emph{every} unfilled fort.  It is possible that unfilled components can be involved in multiple forts.  In other words, as long as there is an unfilled fort the oracle can stop the player from filling all vertices.

Proposition \ref{prop:filler} describes a useful property for the $Z_1$-Game on a tree.

\begin{proposition}\label{prop:filler}
Consider the $Z_1$-Game on a tree $T$.  If any vertices are filled, then at no cost the player can get to a point where either all vertices are filled or there is exactly one filled vertex adjacent to two or more unfilled vertices.
\end{proposition}
\begin{proof}
If the conclusion does not hold, then either the filling rule can be applied or there exists two unfilled disjoint subtrees which are not adjacent to a common filled vertex.  In the latter case, hand these components to the oracle, whatever is returned can have at least one vertex filled.  Continue until no more filling can occur and at this point the conclusion holds.
\end{proof}

We now have the tools to describe which trees have $Z_1(T)=2$.

\begin{definition}
A \emph{comb graph} is a tree with maximum degree three and all degree-three vertices are on a single path.  A pair of \emph{initial vertices} of a comb graph are two vertices $u,v$ such that every degree-three vertex is an internal vertex on the unique path between $u$ and $v$.
\end{definition}

\begin{figure}[ht!]
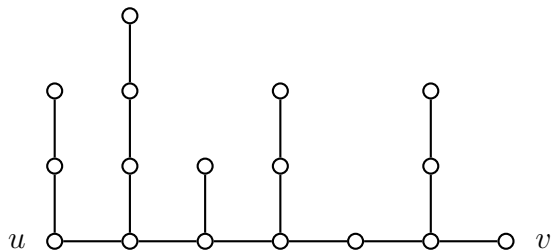

\centering
\COMB
\caption{A comb graph and a pair of initial vertices $u$ and $v$.}
\label{fig:comb}
\end{figure}

\begin{lemma} \label{lem:comb}
Let $G$ be a tree with $Z_1(G)=2$.  Then it is a comb graph.  Moreover, the two tokens must be spent on a pair of initial vertices.
\end{lemma}
\begin{proof}
Since $Z_1(G)>1$,  Theorem~\ref{thm:Z0} says that $G$ is not a path so $G$ must have a vertex with degree at least three.

Lemma \ref{lem:spendk} says that if only one token is spent then any use of the oracle can result in a situation where no progress is made, thus only the filling rule can be applied or the second token can be spent.  However, the operations of spending the second token and applying the filling rule can be interchanged so it can be assumed, without loss of generality, that the two tokens are spent initially.  (In contrast to the previous example, here we can assume that we spend both up front as we do not use the oracle to make progress.)


If the graph has a vertex $x$ of degree at least four, then after spending the two tokens there are at least two disjoint subtrees adjacent to $x$ which contain no filled vertices, these subtrees form a fort.

If not all degree three vertices are on the unique path joining the two initial filled vertices, then there is a degree three vertex which contains two disjoint subtrees which contain no filled vertices, these disjoint subtrees also form a fort.

Combining the previous two observations it can be concluded that $G$ must be a comb graph and that the two places where tokens were spent are initial vertices.

It remains to show that for any comb graph if tokens are spent on any pair of initial vertices the entire graph can be filled. For this, appeal to Proposition~\ref{prop:filler}.  In particular, if not all vertices are filled, then there would have to be a vertex neighboring two or more unfilled disjoint subtrees.  But this implies either the existence of a vertex of degree at least four or a vertex of degree at least three not on the path joining the initial vertices; both of which are impossible for the comb graph.
\end{proof}

A \emph{zig-zag path} is an outerplanar graph which is not a path and can be decomposed as a pair of paths with at least one additional edge between the paths (for a precise definition, see \cite{PropTime}). If spending at two vertices on a zig-zag path make the whole graph filled in the $Z$-Game, then the corresponding \emph{ending vertices} are the vertices which never have the filling rule used to fill a neighbor.

\begin{definition}
A \emph{pick comb graph} is a graph obtained from either combining a zig-zag path and a comb by identifying the ending vertices of the zig-zag path to a pair of initial vertices of the comb, and/or by connecting the ends of a path to a pair of initial vertices of the comb. 
\end{definition}

\begin{figure}[ht!]
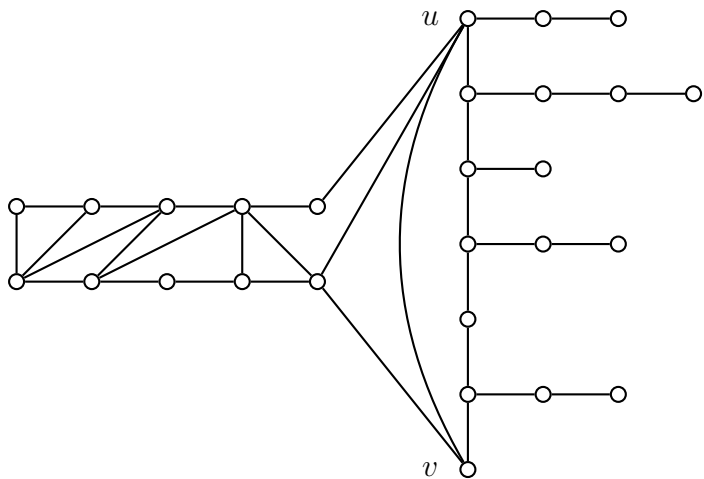

\begin{center}
\PICKCOMB
\end{center}
\caption{A pick comb graph.}
\label{fig:pickcomb}
\end{figure}

\begin{theorem}
Let $G$ be a connected graph.  Then we have the following.
\begin{itemize}
\item $Z_0(G)=2$ if and only if exactly one block of $G$ has a cycle and $G$ does not have $K_4$ or $T_3$ as a minor.
\item $Z_1(G)=2$ if and only if $G$ is a zig-zag path, comb graph, or pick comb graph.
\item $Z_q(G)=2$ for $q\ge 2$ if and only if $G$ is a zig-zag path.
\end{itemize}
\end{theorem}

\begin{proof}
The result for $Z_0(G)$ follows from \cite{EGR}.  When $q\ge 2$, $Z_q(G)=2$ if and only if $Z(G)=2$ by Theorem~\ref{thm:spendfirst}.  Then, $Z(G)=2$ if and only if $G$ is a zig-zag path by \cite{PropTime,Row,Yang13}.   All that remains is to determine what happens for $Z_1(G)$.

By the definition of a pick comb graph, a zero forcing set can be found for the zig-zag path (or just path) that, upon applying the filling rule, ends at initial vertices of the comb graph.  Using the oracle and the filling rule the remaining vertices can be filled as discussed in Lemma \ref{lem:comb}.  In particular, all pick comb graphs have $Z_1(G)\le 2$.  Moreover, if $Z_1(G)=1$, then $Z(G)=1$ and the graph would have to be a path. Thus, $Z_1(G)=2$ for all pick comb graphs.

Conversely, suppose $G$ is a graph with $Z_1(G)=2$.  As in the proof of Lemma~\ref{lem:comb} assume the first two actions of the player are to spend tokens and then apply the filling rule until no more vertices can be filled.  If all vertices have been filled, then the graph is a zig-zag path.  Now assume that there remain unfilled vertices.

At this point, let $u$ and $v$ be the two filled vertices which were not used by the filling rule to fill another vertex.  Let $Y$ be the induced subgraph on the set of filled vertices.  Then $Y$ consists of several possibilities:
\begin{itemize}
\item A pair of paths, possibly consisting of one vertex.
\item A single path.
\item A zig-zag path.
\end{itemize}
Let $H$ be the graph obtained from $G$ by removing all filled vertices except for $u$ and $v$ and removing the edge $\{u,v\}$, if it exists.  Since every vertex that performed a force is not adjacent to any current unfilled vertex, $\{u,v\}$ is a cut-set of $G$.  Therefore, $G=Y\cup H$ and $V(Y)\cap V(H)=\{u,v\}$.

Now consider the structure of $H$.  We can use $\{u,v\}$ as a cut-set to partition the graph; for this we look at the set of components that are in $H-\{u,v\}$ where we then add back copies of $u$ and $v$ in each resulting component where an adjacency occurred.  Define $L$ to be the components of $H$ that contain $u$ but not $v$; $R$ to be the components of $H$ that contains $v$ but not $u$; and $M$ to be the components of $H$ that contain both $u$ and $v$.  Each component in $L$ and $R$ has $Z_1$ equal to $1$, so is a path.  If any of $|L|$, $|R|$, or $|M|$ are at least two, then these components can be used to form a fort which contradicts $Z_1(G)=2$.

Therefore, $|L|$, $|R|$, and $|M|$ are each at most one. If $|M|=0$, then at least one of $L$ and $R$ must consist of at least two paths since there remain unfilled vertices so at least one of $u$ or $v$ is adjacent to two or more unfilled vertices.  But these correspond to a pair of unfilled paths sharing a common filled vertex as a neighbor which forms a fort.  This is impossible.

Using a similar argument it can be concluded that $|M|=|L|=|R|=1$.  For the remainder of the proof let $M$, $L$, and $R$ denote these corresponding, single components.  

If $u$ has at least two neighbors in $M$, then $V(M)\cup V(R)\setminus\{u,v\}$ is a fort, so $u$ is a leaf in $M$.  Similarly, $v$ is also a leaf in $M$.  Let $u_1$ and $u_2$ be the neighbor of $u$ in $L$ and in $M$, respectively.  Whenever the player gives two components to the oracle, the oracle can either return a component that does not contain $u_1$ and $u_2$, or return two components that contain $u_1$ and $u_2$.  In this way the vertex $u$ cannot perform any force until $u_2$ is filled.  Let $H'$ be the graph obtained from $H$ by removing $u$ and the branch $L$.  This means the vertex $u$ has no effect on any force that happens in $H'$, and $Z_0(H')=1$ by using $\{v\}$ as the zero forcing set.  Therefore, $H'$ is a tree.  Since $u$ is a leaf on $M$ and $L$ is a path, the graph $H$ is a tree.

Since $H$ is a tree with $Z_1(H)=2$ then, by Lemma~\ref{lem:comb}, $H$ is a comb graph with $u$ and $v$ a pair of initial vertices of $H$.  Thus, $G$ is a pick comb graph as desired.
\end{proof}

\section{Efficient computation of \texorpdfstring{$Z_1(T)$}{Z1(T)} for \texorpdfstring{$T$}{T} a tree}\label{sec:Z1onTree}

Due to the nature of the $Z_q$-Game, it is not immediately clear how to go about computing $Z_q(G)$ for a given graph.  An algorithm was presented in the paper introducing the $Z_q$-Game (see \cite{BGH}), which runs over all $2^{|V(G)|}$ subsets of the graph, only allowing for computation on graphs of moderate size (see \cite{BGH}). If $S$ is an induced subgraph of $G$, then for any vertex $v$ in $S$, we let $\deg_S(v)$ denote the degree of $v$ in the induced subgraph $S$. 

We now give a method to compute $Z_1(T)$ when $T$ is a tree which can be determined in polynomial time.  We start with the following characterization for $Z_1(T)$.

\begin{theorem}\label{thm:Z1onTree}
Let $T$ be a tree with $|T|\ge 3$. Let $\mathcal{P}_S(v)$ denote the set of maximal paths in a tree $S$ with $v$ as one endpoint (i.e., paths with $v$ as an endpoint and not contained in a longer path).  Then
\begin{equation}\label{eq:Z1}
Z_1(T)=2+\max_{v\in V(T)}\bigg(\max_{\substack{S\subseteq T\\ v\in V(S)}}\bigg(\min_{P\in\mathcal{P}_S(v)}\bigg(\sum_{w\in V(P)}(\deg_S(w)-2)\bigg)\bigg)\bigg).
\end{equation}
\end{theorem}
The proof will be split into two parts, namely establishing the right hand side as lower and upper bound for $Z_1(T)$.

\subsection{Lower bound for $Z_1(T)$}

Fix a vertex $v$ and a corresponding subtree $S$ achieving  the maximums in \eqref{eq:Z1}.

We modify the game by making the oracle more generous in the following ways:
\begin{itemize}
\item The oracle announces the vertex $v$ and the subtree $S$ to the player. 
\item The oracle announces that if any vertex in $S$ is filled which is adjacent to a vertex in $T\setminus S$, then the corresponding vertex in $T\setminus S$ is filled as well as the \emph{entire} corresponding subtree in $T\setminus S$ containing that vertex.
\item The oracle fills in vertex $v$ before the game begins (and any corresponding subtrees in $T\setminus S$ to which $v$ is adjacent).
\item Whenever a token is spent the oracle will fill in all vertices in the unfilled subtree where the token was spent off of the unique filled vertex which is adjacent to that subtree.  And then apply the filling rule if possible.
\end{itemize}
Note that all of these new options only offer an advantage to the player in filling the vertices.  So any lower bound in this modified version of the game is also a lower bound in the original version.

Using the information above, we may assume the player does not spend a token in $T\setminus S$ as a better effect can be reached by spending in the vertex in $S$ next to the corresponding subtree in $T\setminus S$.

Label the initial vertex $v$ as $v_0$, and the game play proceeds as follows with the action of the player and responses from the oracle:
\begin{itemize}
\item In the first $\deg_S(v_0)-1$ rounds the player will fill a vertex and the oracle will then fill in all vertices of the subtree with root at $v_0$ containing the newly filled vertex, along with any adjacent parts of $T\setminus S$.  

After $\deg_S(v_0)-1$ rounds,  there is precisely one unfilled neighbor of $v_0$, denote this neighbor as $v_1$, which, by the filling rule, can now be filled.  Now consider vertex $v_1$.
\item Assume the process has taken us to vertex $v_i$ with $i\ge1$.  The player will spend tokens and the oracle will fill in all of the subtrees rooted at $v_i$ in which the token was spent, along with any adjacent parts of $T\setminus S$.

This continues until $v_i$ has one unfilled neighbor (this occurs after spending $\deg_S(v_i)-2$ tokens), at which point the filling rule is applied and vertex $v_{i+1}$ will be considered.
\item Once a vertex $v_k$ in $S$ is reached which is a leaf, all vertices have been filled and the game ends.
\end{itemize}

(In the above analysis the reason $\deg_S(v_i)$ is used and not $\deg_T(v_i)$ is because anything which is in $T\setminus S$ is automatically filled.)

Due to the nature of the play of the oracle, the only option the player can ever use is spending a token (e.g., all unfilled components always share a common vertex and all applications of the filling rule are done by the oracle).  Moreover, the above places no restrictions on where tokens are spent.


We now determine the cost to the player to win the altered game using the above strategy.  The sequence of vertices $v_i$ forms a maximal path starting at $v_0$ and going to $v_k$ in $S$.  Moreover, the cost is given by
\[
(\deg_S(v_0)-1)+(\deg_S(v_1)-2)+\cdots+(\deg_S(v_{k-1})-2)=2+\sum_{i=0}^k(\deg_S(v_i)-2)
\]
%
%
(the $2$ in front comes from a combination of correcting the first term as well as accounting for the leaf on the end of the path).

Finally, note that the player can choose the path in $S$ by spending tokens in the subtrees \emph{off} the path.  So the minimum number of tokens needed for the player to win is thus
\[
2+\min_{P\in\mathcal{P}_S(v)}\bigg(\sum_{w\in V(P)}(\deg_S(w)-2)\bigg).
\]
The lower bound now follows.

\subsection{Upper bound for $Z_1(T)$}

\subsubsection{An important algorithm}

The upper bound discussion begins with a recursive Algorithm~\ref{algorithm} for which a tree $T$ rooted at $v$ assigns a value, $f_{T,v}(w)$, to each vertex $w$.

\bigskip

\begin{algorithm}[H]
\SetAlgoLined
\eIf{$w$ has no descendents}{
$f_{T,v}(w):=0$
}{
let $u_0,\ldots,u_k$ be the children of $w$\;
let $T_i$ be the subtree rooted at $u_i$ with all descendents of $u_i$\;
compute $f_{T_i,u_i}(u_i)$ for all $i$\;
relabel (if needed) so $f_{T_0,u_0}(u_0)\ge f_{T_1,u_1}(u_1) \ge \cdots\ge f_{T_k,u_k}(u_k)$\;
$f_{T,v}(w):=\max\{0+f_{T_0,u_0}(u_0),1+f_{T_1,u_1}(u_1),\ldots,k+f_{T_k,u_k}(u_k)\}$
}
\caption{Computation of $f_{T,v}(w)$ for all vertices $w$ in tree $T$ rooted at $v$}
\label{algorithm}
\end{algorithm}

\bigskip

The value $f_{T,v}(w)$ can be interpreted in the following way:  In the $Z_1$-Game on the tree $T$ rooted at $v$, given that the only unfilled vertices are those which are descendents of the current vertex ($w$), then $f_{T,v}(w)$ is the minimal number of tokens needed to be spent by the player to fill in the remaining vertices.  The justification for this interpretation comes from combining the following two observations.

\begin{observation}
If a vertex is a leaf, then no additional tokens are needed to be spent (i.e., everything has been filled).
\end{observation}

\begin{observation}\label{obs:trim}
If a vertex is not a leaf, and the minimum number of tokens for all subtrees rooted at the children of the vertex have been determined, then the minimum number of tokens for that vertex can be determined.  This is done by arranging the values for the number of tokens of the subtrees of the children in weakly decreasing order, e.g., $c_0\ge c_1\ge\cdots\ge c_k$.  Then, the number of tokens for the vertex is $\max_{0\le i\le k}\{i+c_i\}$.
\end{observation}

The reason this works is because as the player is at a vertex the oracle will allow them to move down a level when it is more expensive to work in a particular subtree than to work in any other available subtree.  Filling a branch will cost one token (by interaction with the oracle, see Proposition~\ref{prop:filler}).  So in the worst case scenario, the oracle will make the player spend $i$ tokens on the most expensive subtrees before going down, where $i$ is chosen to maximize $i+c_i$.

\begin{figure}[ht!]
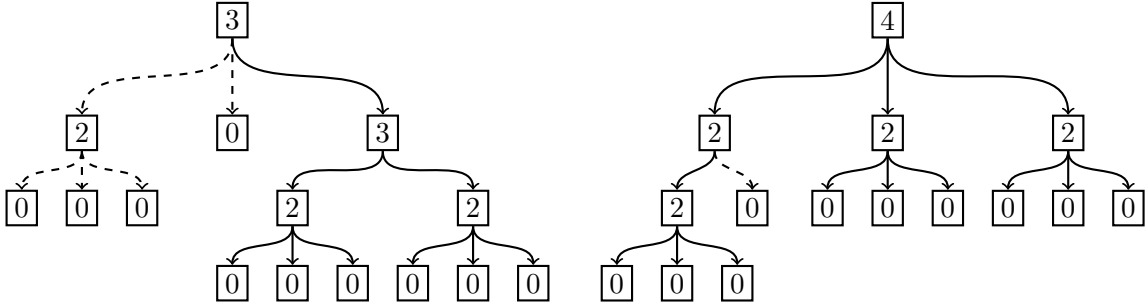

\begin{center}
\ALGA
\hfil
\ALGB
\end{center}
\caption{Two examples of the assignment of minimum cost to vertices assigned to a tree which differ by location of root.}
\label{fig:alg}
\end{figure}

Two examples of the values for the same tree are given in Figure~\ref{fig:alg} (the difference being the location of the roots). The dashed subtrees are those which will not be needed when determining optimal cost at a vertex.

\subsubsection{Properties of the algorithm}
The algorithm has several important properties, that when combined will establish the upper bound for Theorem~\ref{thm:Z1onTree}.  Before proceeding, a simple computation shows for $T=P_n$ that $f_{T,v}(v)=0$ for the leaves and $f_{T,v}(v)=1$ for internal vertices.

\begin{proposition}
For a tree $T$ on at least three vertices we have
\[
f_{T,v}(v) = 2+\max_{\substack{S\subseteq T\\ v\in V(S)}}\bigg(\min_{P\in\mathcal{P}_S(v)}\bigg(\sum_{w\in V(P)}(\deg_S(w)-2)\bigg)\bigg).
\]
\end{proposition}
\begin{proof}
Using the language in Observation~\ref{obs:trim}, if at any vertex $w$ we have $\max\{i+c_i\}$ happen for $i<k$, then prune off the subtrees corresponding to $i+1,\ldots,k$. (As an example, in Figure~\ref{fig:alg} the parts to be pruned are shown with dashed lines).   The resulting tree will be $S$.

Now, starting at the root $v$, to move down to a leaf  $\deg_S(w)-1$ tokens must be spent at the ends of the path and $\deg_S(w)-2$ spent on the interior of the path (if the ``$+2$'' is incorporated then the costs can be $\deg_S(w)-2$ for all vertices on the path).  The value $\deg_S(w)-2$ corresponds to the ``$i$'' term present in ``$i+c_i$'', in particular, the cost from the algorithm of a vertex is found by minimizing the sum $\sum(\deg_S(w)-2)$ over all paths from that vertex down to a leaf which is the right hand side of Equation \ref{eq:Z1}.
\end{proof}

\begin{proposition}
For a tree $T$ on at least three vertices $f_{T,v}(v)+1\ge \max_{w}f_{T,w}(w)$.
\end{proposition}
\begin{proof}
Observe that if $u$ is a child of $x$, then $f_{T,v}(x)\ge f_{T,v}(u)$ (e.g., $f_{T,v}(u)$ is at least the cost of any subtree off the vertex $x$).  From this it can be seen that $f_{T,v}(x)$ is maximized when $x=v$.

Fix a $w$ in $T$ which maximizes $f_{T,w}(w)$ and so that in the corresponding maximal $S$ there are at least two children of $w$.  Such a $w$ can be found by starting at an arbitrary $w$ which maximizes $f_{T,w}(w)$; then find a vertex $x$ which is as far down in the tree as possible and has $f_{T,w}(w)=f_{T,w}(x)$ (possibly $w=x$); $x$ is the desired vertex.

Finally, observe that $f_{T,v}(w)$ is found by looking over the subtrees of $w$ in $T$ rooted at $v$.  But this has exactly one less subtree than the computation when determining $f_{T,w}(w)$ and the removal of one subtree can only decrease the value by at most one, e.g., $f_{T,v}(w)\ge f_{T,w}(w)-1$.

Combining these ideas yields $f_{T,v}(v)\ge f_{T,v}(w)\ge f_{T,w}(w)-1$, establishing the result.
\end{proof}

\begin{proposition}
If $T$ is a tree on three or more vertices, then $f_{T,v}(v)$ is not constant on the vertices of $T$.
\end{proposition}
\begin{proof}{}

If $T$ is a path on $n$ vertices with $n\geq 3$, then $f_{T,v}(v)=0$ if $v$ is a leaf, otherwise $f_{T,v}(v)=1$.  Thus, $f_{T,v}(v)$ is not constant, so we may assume $T$ has at least a vertex of degree three or more.

Define $g_u(v)=f_{T^u_v,v}(v)$, where $T^u_v$ is the component of $T-\{u,v\}$ that contains vertex $v$. Then define a digraph $\Gamma$ on the vertex set $V(T)$.  For each edge $uv\in E(T)$, add an arc $(u,v)$ to $\Gamma$ if $g_v(u)\leq g_u(v)$, and add an arc $(v,u)$ to $\Gamma$ if $g_u(v)\leq g_v(u)$, so $\Gamma$ is doubly directed on vertices $\{u,v\}$ if $g_v(u)=g_u(v)$.

Pick a strongly connected component of $\Gamma$ that has no out neighbor, and let $S$ be the induced subgraph of $T$ on this strongly connected component.  Let $u_0$ be a leaf of $S$ or the unique vertex of $S$ if $|V(S)|=1$.  Since $T$ has a vertex of degree at least $3$, for each leaf $u$ and its unique neighbor $v$, $g_u(v)\geq 1$ and $g_v(u)=0$.  Thus, $u_0$ is not a leaf of $T$, and $N_T(u_0)\setminus N_S(u_0)$ is not empty.  Let $M=f_{T,u_0}(u_0)$.  We will find a vertex $v\in N_T(u_0)\setminus N_S(u_0)$ such that $f_{T,v}(v)\leq M-1$.

If $v\in N_S(u_0)$, then $g_v(u_0)=g_{u_0}(v)$ by the definition of $\Gamma$, and 
\[g_{v}(u_0)=\max\{0+g_{u_0}(w_0),\ldots,k+g_{u_0}(w_k)\},\]
where $w_0,\ldots,w_k$ are the neighbors of $u_0$ in $T^v_{u_0}$ with $g_{u_0}(w_0)\geq \cdots \geq g_{u_0}(w_k)$.  Thus, $g_{u_0}(v)\geq g_{u_0}(w_0)$, so 
\[\begin{aligned}
f_{T,u_0}(u_0) &= \max\{0+g_{u_0}(v),1+g_{u_0}(w_0),\ldots, (k+1)+g_{u_0}(w_k)\}\\
 &= \max\{g_{u_0}(v),1+g_{v}(u_0)\}=1+g_v(u_0).
 \end{aligned}\]
 Therefore, $g_{u_0}(v)=g_v(u_0)=M-1$.  On the other side, if $v\in N_T(u_0)\setminus N_S(u_0)$, then 
 \[g_{u_0}(v)<g_v(u_0)\leq f_{T,u_0}(u_0)=M,\]
so $g_{u_0}(v)\leq M-1$.  In summary, $g_{u_0}(v)\leq M-1$ for all $v\in N_T(u_0)$.

Now let $v_0,\ldots, v_d$ be the neighbors of $u_0$ in $T$ such that $g_{u_0}(v_0)\geq\cdots\geq g_{u_0}(v_d)$.  Then 
\[f_{T,u_0}(u_0) = \max\{0+g_{u_0}(v_0),\ldots, d+g_{u_0}(v_d)\}.\]
Let $\ell$ be the smallest index such that $M=\ell+g_{u_0}(v_\ell)$.  Then 
\[g_{v_\ell}(u_0) = \max\{0+g_{u_0}(v_0),\ldots, (\ell-1)+g_{u_0}(v_{\ell-1}), \ell+g_{u_0}(v_{\ell+1}), \ldots, (d-1)+g_{u_0}(v_d)\}\leq M-1.\]

Recall that $g_{u_0}(v)\leq M-1$ for all $v\in N_T(u_0)$, so $\ell\geq 1$.  Since $u_0$ is a leaf in $S$, $v_\ell$ is in $N_T(u_0)\setminus N_S(u_0)$.  It follows that $g_{u_0}(v_\ell)<g_{v_\ell}(u_0)$.  Let $w_0,\ldots,w_k$ be the neighbors of $v_\ell$ in $T^{u_0}_{v_\ell}$ such that $g_{v_\ell}(w_0)\geq\cdots\geq g_{v_\ell}(w_k)$.  Then 
\[g_{u_0}(v_\ell)=\max\{0+g_{v_\ell}(w_0), \ldots, k+g_{v_\ell}(w_k)\}\]
and 
\[\begin{aligned}
f_{T,v_\ell}(v_\ell) &= \max\{0+g_{v\ell}(u_0),1+g_{v_\ell}(w_0), \ldots, (k+1)+g_{v_\ell}(w_k)\} \\
 &= \max\{g_{v_\ell}(u_0),1+g_{u_0}(v_\ell)\}=g_{v_\ell}(u_0)\leq M-1.
 \end{aligned}\]
This completes the proof.
 \end{proof}

\begin{proposition}
There is a strategy for the player to fill all vertices in a tree $T$ on at least three vertices using at most $1+\min_v f_{T,v}(v)=\max_v f_{T,v}(v)$ tokens.
\end{proposition}
\begin{proof}
That $1+\min_v f_{T,v}(v)=\max_v f_{T,v}(v)$ follows by combining the two preceding propositions.  So we now only need to give a strategy for the player.  The player spends the first token at a vertex $v$ which achieves $\min_v f_{T,v}(v)$ (this is the ``$+1$'') and computes the values $f_{T,v}(w)$ for all vertices $w$ in the tree.  Every subsequent token is spent by the following rule:
\begin{quote}
From the unique current filled vertex with unfilled children spend in a leaf which is reached by going from the current vertex and always choosing a vertex which is maximal for $f_{T,v}(w)$ among the unfilled children.
\end{quote}
After spending each token the player then uses Proposition~\ref{prop:filler}.  This will result in having the unique filled vertex with unfilled children ending up somewhere on that path where the most ``costly'' subtree not previously filled has now been completely filled.

We analyze the response of the oracle (e.g., how the oracle will respond when using Proposition~\ref{prop:filler}).  We can simplify the decision making by putting it into one of two options:
\begin{itemize}
\item Fill in the entire subtree where the token was spent (and stop).
\item Fill in the subtrees where the token was not spent, go down one level, and then restart the decision making process.
\end{itemize}

The oracle will choose the first option only when the cost to the player will be higher if forced to go down into a different subtree.  This is determined by looking at the costs of filling in the subtrees ($f_{T,v}(w)$) combined with the cost of filling in all of the subtrees which will not be used.  This is precisely the $i+c_i$ expression coming from the algorithm.  

In particular, the total cost after spending the first token will be $f_{T,v}(v)$.
\end{proof}

\subsubsection{Running the algorithm}
When carrying out the algorithm for a specific vertex $v$ we can work from the leaves to the root and never need to visit any vertex more than one time, so is polynomial in the number of vertices.  Since we run the algorithm over all possible vertices we can then conclude that the runtime to determine $Z_1(T)$ is polynomial in the number of vertices of the tree.  

This allows us to compute $Z_1(T)$ for all small trees.  In Table~\ref{tab:Z1data}  the columns correspond to a value of $Z_1(T)$, the rows correspond to the number of vertices, and an entry is how many trees on that number of vertices has that value for $Z_1(T)$.

\begin{table}[!ht]
\[
\begin{array}{l|ccccccccccc}
&k{=}1&k{=}2&k{=}3&k{=}4&k{=}5&k{=}6&k{=}7&k{=}8&k{=}9&k{=}10&k{=}11\\ \hline 
n{=} 3& 1 &&&&&&&&&&\\ 
n{=} 4& 1&1 &&&&&&&&&\\ 
n{=} 5& 1&1&1 &&&&&&&&\\ 
n{=} 6& 1&3&1&1 &&&&&&&\\ 
n{=} 7& 1&5&3&1&1 &&&&&&\\ 
n{=} 8& 1&10&7&3&1&1 &&&&&\\ 
n{=} 9& 1&17&17&7&3&1&1 &&&&\\ 
n{=}10& 1&35&39&19&7&3&1&1 &&&\\ 
n{=}11& 1&63&95&45&19&7&3&1&1 &&\\ 
n{=}12& 1&126&228&118&47&19&7&3&1&1 &\\ 
n{=}13& 1&240&559&298&125&47&19&7&3&1&1 \\ 
n{=}14& 1&479&1372&781&321&127&47&19&7&3&1\\ 
n{=}15& 1&934&3387&2031&855&328&127&47&19&7&3\\ 
n{=}16& 1&1867&8399&5372&2266&880&330&127&47&19&7\\ 
n{=}17& 1&3687&20871&14223&6081&2344&887&330&127&47&19\\
n{=}18& 1&7372&52010&38002&16353&6336&2369&889&330&127&47\\
n{=}19&1& 14654& 129792& 101844& 44312& 17136& 6416& 2376& 889& 330& 127\\
n{=}20& 1& 29304& 324514& 274449& 120437& 46721& 17396& 6441& 2378& 889& 330
\end{array}
\]
\caption{The number of trees on $n$ vertices with $Z_1(T)=k$.}
\label{tab:Z1data}
\end{table}

\subsubsection{Another view of computing $Z_{1}(T)$}
As noted in Theorem \ref{thm:Z1onTree} we have developed a recursive formulation for computing $Z_{1}(T)$ for a given tree $T$. As a matter of completeness, we offer a second viewpoint for computing $Z_{1}(T)$, which is also recursive in nature, but may be of general interest.  

 Suppose $T$ is a tree that is not a path. Then let
 $L_T$ be the set of leaves of $T$, and let $u_1$ and $u_2$ be two distinct vertices in $L_T$. 
Let $U$ be the set of internal vertices of the 
unique path in $T$ from $u_1$ to $u_2$. Suppose $v \in U$ and denote by $S$
the tree obtained from $T$ by deleting 
the vertices of the two components in $T-v$ that contain 
$u_1$ or $u_2$, respectively. As before we let $\mathcal{P}_S(v)$ denote the set of 
maximal paths in the tree $S$ with $v$ as one endpoint. Then
\[
Z_1(T) = 2 + \min_{u_1,u_2 \in L_T} (\max_{v \in U} (\min_{P \in \mathcal{P}_S(v)} (\sum_{w \in V(P)} (\deg_S(w)-2)))). \]
We omit the verification of this equation, as it resembles the proof that we used for the previous interpretation. We also note here that it can be shown that for a given tree $T$ on $n$ vertices $Z_{1}(T)$ can be computed in $O(n^5)$ time.

\section{Concluding remarks}\label{sec:conclusion}

Our main interest has been on the combinatorial aspects of the $Z_q$-Game.  Given the nature of the need to find consistent strategies there is still much that is not known about this game, including determining $Z_q(G)$ for many well-known families of graphs.  In addition, we list several problems here.
\begin{itemize}
\item As noted in the introduction the $Z$-Game works by looking at sets, and not strategies, since all spending of tokens can occur at the start.  There are some graphs for which all spending in the $Z_q$-Game can occur at the start, e.g., any graph where $Z_q(G)\le q$.  Are there any family of graphs where $Z_q(G)>q$ and we can still spend all tokens at the start?
\item We found an efficient way to compute $Z_1(T)$ for a tree $T$ without having to directly run the game.  Are there similar methods that work for other general graphs?  (The determination of $Z(G)$ is believed to be NP, so the answer would be expected to be no in general, but perhaps for some special subfamilies progress can be made.)
\item For trees $T$ we now have ways to compute $Z_0(T)$ (always $1$), $Z_1(T)$ (see Theorem~\ref{thm:Z1onTree}), and $Z_\infty(T)$ (path cover number).  Do there exist other ``simple'' combinatorial interpretations that can be used to compute the $q$-analogue zero forcing number for trees for other values of $q$?
\item The original motivation for the study of the $Z_q$-Game comes from its connection to bound nullity of a matrix associated with a graph (and more generally the inertia).  We can see even for trees that the parameter associated with the $Z_1$-Game can be arbitrarily far off from the nullity for a graph.  As an example, a complete binary tree of depth $d$ has $Z_1(T)=d$,  but any matrix associated with the graph has nullity at most $2$.  Are there additional ways to modify the game to get improved bounds on nullity?
\end{itemize}

\section*{Acknowledgments:}
The genesis of this paper was from the American Institute of Mathematics workshop \emph{Zero forcing and its applications}, held in January 2017.  We thank AIM for hosting their workshop and their support in bringing us together.  

Steve Butler was partially supported by a grant from the Simons Foundation (\#427264).

Shaun Fallat's research is supported in part by an NSERC Discovery Research Grant, Application No.: RGPIN-2014-06036.

\end{document}